\documentclass[11pt]{article}
\usepackage{amsthm, amsfonts, amsmath}
\usepackage{amscd}
\usepackage[top=27truemm,bottom=35truemm,left=20truemm,right=20truemm]{geometry}

\allowdisplaybreaks[4]
\usepackage{graphicx}

\makeatletter
\@addtoreset{equation}{section}
\def\theequation{\thesection.\arabic{equation}}
\makeatother

\def\g{{\rm gcd}}
\def\pp{\prod_{l\in{\cal P}(q)}\left(1-\frac{1}{l^2}\right)}
\def\p{\hspace{-2pt}+\hspace{-2pt}}
\def\m{\hspace{-1pt}-\hspace{-1pt}}
\def\tables{\begin{table}[!t]
\begin{center}
\begin{tabular}{c||c|rccc|lc}
$\frac{x}{z}$&$n$&$xw-yz=1$&example of $w^2$&$\g(9,z)$&$k$&\hspace{11mm}conditions of $m$&$m$\\[1mm] \hline \hline
$\frac{1}{0}$&1&$w=1$&1&9&1&$1=\g(9,m)$ and $m\equiv1\,\mbox{mod}\,9$&1\\[1mm]
$\frac{2}{9}$&1&$2w-9z=1$&25&9&7&$1=\g(9,m)$ and $7m\equiv1\,\mbox{mod}\,9$&4\\[1mm]
$\frac{4}{9}$&1&$4w-9z=1$&4&9&4&$1=\g(9,m)$ and $4m\equiv1\,\mbox{mod}\,9$&7\\[1mm] \hline
$\frac{1}{3}$&3&$w-3y=1$&1&3&1&$3=\g(9,m)$ and $\frac{m}{3}\equiv1\,\mbox{mod}\,3$&3\\[1mm]
$\frac{2}{3}$&3&$2w-3y=1$&1&3&1&$3=\g(9,m)$ and $\frac{m}{3}\equiv1\,\mbox{mod}\,3$&3
\end{tabular}
\caption{Rotation numbers and exponents about $X_9$}
\end{center}
\end{table}
\begin{table}[!t]
\begin{center}
\begin{tabular}{c||c|rccc|lc}
$\frac{x}{z}$&$n$&$xw-yz=1$&example of $w^2$&$\g(10,z)$&$k$&\hspace{11mm}conditions of $m$&$m$\\[1mm] \hline \hline
$\frac{1}{0}$&1&$w=1$&1&10&1&$1=\g(10,m)$ and $m\equiv1\,\mbox{mod}\,10$&1\\[1mm]
$\frac{3}{10}$&1&$3w-10z=1$&9&10&9&$1=\g(10,m)$ and $9m\equiv1\,\mbox{mod}\,10$&9\\[1mm] \hline
$\frac{1}{5}$&2&$w-5y=1$&1&5&1&$2=\g(10,m)$ and $\frac{m}{2}\equiv1\,\mbox{mod}\,5$&2\\[1mm]
$\frac{2}{5}$&2&$2w-5y=1$&4&5&4&$2=\g(10,m)$ and $4\frac{m}{2}\equiv1\,\mbox{mod}\,5$&8\\[1mm] \hline
$\frac{1}{2}$&5&$w-2y=1$&1&2&1&$5=\g(10,m)$ and $\frac{m}{5}\equiv1\,\mbox{mod}\,2$&5\\[1mm]
$\frac{1}{4}$&5&$w-4y=1$&1&2&1&$5=\g(10,m)$ and $\frac{m}{5}\equiv1\,\mbox{mod}\,2$&5
\end{tabular}
\caption{Rotation numbers and exponents about $X_{10}$}
\end{center}
\end{table}
\begin{table}[!t]
\begin{center}
\begin{tabular}{c||c|rccc|lc}
$\frac{x}{z}$&$n$&$xw-yz=1$&example of $w^2$&$\g(12,z)$&$k$&\hspace{11mm}conditions of $m$&$m$\\[1mm] \hline \hline
$\frac{1}{0}$&1&$w=1$&1&12&1&$1=\g(12,m)$ and $m\equiv1\,\mbox{mod}\,12$&1\\[1mm]
$\frac{5}{12}$&1&$5w-12z=1$&25&12&1&$1=\g(12,m)$ and $m\equiv1\,\mbox{mod}\,12$&1\\[1mm] \hline
$\frac{1}{6}$&2&$w-6z=1$&1&6&1&$2=\g(12,m)$ and $\frac{m}{2}\equiv1\,\mbox{mod}\,6$&2\\[1mm] \hline
$\frac{1}{4}$&3&$w-4y=1$&1&4&1&$3=\g(12,m)$ and $\frac{m}{3}\equiv1\,\mbox{mod}\,4$&3\\[1mm]
$\frac{3}{4}$&3&$3w-4y=1$&1&4&1&$3=\g(12,m)$ and $\frac{m}{3}\equiv1\,\mbox{mod}\,4$&3\\[1mm] \hline
$\frac{1}{3}$&4&$w-3y=1$&1&3&1&$4=\g(12,m)$ and $\frac{m}{4}\equiv1\,\mbox{mod}\,3$&4\\[1mm]
$\frac{2}{3}$&4&$2w-3y=1$&1&3&1&$4=\g(12,m)$ and $\frac{m}{4}\equiv1\,\mbox{mod}\,3$&4\\[1mm] \hline
$\frac{1}{2}$&6&$w-2y=1$&1&2&1&$6=\g(12,m)$ and $\frac{m}{6}\equiv1\,\mbox{mod}\,2$&6
\end{tabular}
\caption{Rotation numbers and exponents about $X_{12}$}
\end{center}
\end{table}
}

\newtheorem{theo}{Theorem}[section]

\newtheorem{lem}[theo]{Lemma}
\newtheorem{rem}[theo]{Remark}
\newtheorem{prop}[theo]{Proposition}
\newtheorem{cor}[theo]{Corollary}
\newtheorem{nota}[theo]{Notation}
\newtheorem{ex}[theo]{Example}

\title{Equations associated to the Principal\\Congruence Subgroup of Level Eight}
\author{Takashi Niwa\footnote{Department of Mathematics, Tokyo Institute of Technology, 2-12-1 Oh-okayama, Meguro-ku, Tokyo 152-8551, Japan, e-mail: niwa.t.ad@m.titech.ac.jp}}
\date{}
\begin{document}
\maketitle

\section*{Introduction}
For a natural number $q$, let $\Gamma_q$ be the principal congruence subgroup of level $q$ of $SL\left(2,\mathbb{Z}\right)$, namely,
\[\Gamma_q=\left\{\begin{pmatrix}
a&b\\c&d
\end{pmatrix}
\in SL\left(2,\mathbb{Z}\right)\,|\,
a\equiv d\equiv1,\,b\equiv c\equiv0\,\,({\rm mod}\ q)\right\}.\]
We consider the compactification of $\mathbb{H}/\Gamma_q$ and denote this compact Riemann surfaces by $X_q$.
It turns out that all compact Riemann surfaces are obtained by algebraic equations.
Thus, $X_q$ is also obtained by an algebraic equation.
Equations of $X_q$ for prime numbers $q$ are given in \cite{II}.
Furthermore, Ishida gives equations for all natural numbers $q$ in \cite{Is}.
They consider a family of modular functions
\[X_r(\tau)=\exp\left(2\pi\sqrt{-1}\cdot\frac{(r-1)(q-1)}{4q}\right)\prod_{s=0}^{q-1}\frac{K_{r,s}(\tau)}{K_{1,s}(\tau)}\]
for $r\in\mathbb{Z}$ which aren't divided by $q$.
Here, $K_{u,v}(\tau)$ are Klein forms of level $q$.
For Klein forms, let see Kubert and Lang \cite{KL}.
They show that $X_3(\tau)$ is integral over $\mathbb{Q}\left[X_2(\tau)^{\varepsilon_q}\right]$, where $\varepsilon_q$ is 1 or 2 according to whether $q$ is odd or even and then, get equations of $X_q$.
Recently, Yang gives another way to get the equations in \cite{Ya} by generalizing Dedekind $\eta$-functions.

Then the purpose of this paper is to find an equation associated to $X_8$ without considering modular functions.
We remark that equations for $q\le7$ are obtained by easy ways.
For $q\le5$, equations of $X_q$ is $y=0$ since the genus of $X_q$ is zero so that  $X_q$ is isomorphic to the Riemann sphere.
An equation of $X_6$ is an elliptic curve
\[y^2=x^3-1\]
since $X_6$ has an automorphism of order 3 with fixed points.
The Klein's quartic $X_7$ (cf. \cite{Kl}) is given by the classical equation
\[y^7=x(x-1)^2.\]
The advantage of our way is that it is resolved in a more simpler way and we can find some properties of compact Riemann surfaces from equations.

We give a summary of the contents.

In \S1, we see some properties of compact Riemann surfaces defined by an equation
\[y^p=\prod_{i=1}^r(x-a_i)^{m_i}.\]
We call it a "semi-hyperelliptic curve".
Kuribayashi showed in \cite{Ku} that if a compact Riemann surface $X$ has an automorphism $\tau$ such that the genus of $X/\langle\tau\rangle$ is zero,  
an equation of $X$ is given by a semi-hyperelliptic curve.
Here, $p$ is the order of the automorphism $\tau$ and $r$ is given by the number of branched points of the natural projection $X\rightarrow X/\langle\tau\rangle$.
Girondo and Gonz\^alez-Diez give the values of the exponents $m_i$ for prime number $p$ in \cite{GG} by considering the rotation number of the automorphism $\tau$.
Using their idea, we obtain the values of the exponents $m_i$ for all natural numbers $p$.
We also obtain some conditions of $a_i$.

In \S2, we have the genus of $X_q/\langle z\mapsto z+n\rangle$ and see that there is an automorphism $\tau$ of $X_q$ such that the genus of $X_q/\langle\tau\rangle$ is zero, for $q\le10$ or $q=12$.
Thus, for these $q$, especially	$q=8$, $X_q$ is obtained by semi-hyperelliptic curves.

In \S3, we consider equations of $X_q$.
Then we give a new approach to get equations of $X_q$ for $q\le10$ or $q=12$, which are semi-hyperelliptic curves.
By using a method of \S1, we get equations of these $X_q$ except for values of constant numbers $a_i$.
It is done by computing rotation numbers of parallel displacement.

In \S4, we get an equation of $X_8$ completely including constant numbers $a_i$.
The equation is
\[y^8=x^2(x-1)(x+1).\]
We determine $a_i$ from the automorphism group of $X_8$ agree with the one of the compact Riemann surface defined by an equation which we get in \S3.
We consider automorphisms in the projective space because it is difficult to consider them on algebraic curves in $\mathbb{C}^2$.
\section{Properties of semi-hyperelliptic curves}
We recall that a relation between algebraic functions and compact Riemann surfaces.
Let
\[F(x,y)=\sum_{i=0}^pa_i(x)\cdot y^i\in\mathbb{C}[x,y]\]
be an irreducible polynomial.
If $p\ge1$, there exists a compact Riemann surface $X_F$ which contains the connected Riemann surface
\begin{eqnarray}\label{eq:com}
\left\{\,(x,y)\in\mathbb{C}^2\,|\,F(x,y)=0,F_y(x,y)\ne0,a_p(x)\ne0\,\right\},
\end{eqnarray}
which is an open Riemann surface with finitely many complementary points.
The first projection $(x,y)\mapsto x$ is a holomorphic function and admits a holomorphic extension $X_F\rightarrow \hat{\mathbb{C}}$.
$X_F$ is uniquely determined by $F(x,y)$ up to conformal maps.
Conversely, all compact Riemann surface $X$ has an irreducible polynomial $F(x,y)$ such that $X_F$ is isomorphic to $X$. 
It is shown by considering the meromorphic function field.
Let ${\cal M}(X)$ be the meromorphic function field of $X$.
${\cal M}(X)$ is an algebraic function field of one variable.
In other words, there exists an element $f\in{\cal M}(X)\setminus\mathbb{C}$ such that the extension $\mathbb{C}(f)\subset{\cal M}(X)$ is finite.
It is known that
\[\left[{\cal M}(X):\mathbb{C}(f)\right]={\rm deg}(f)\]
and if $g\in{\cal M}(X)$ is injective on the generic fiber of $f$, we get ${\cal M}(X)=\mathbb{C}(f,g)$.
There is an irreducible polynomial $F(x,y)$ such that $F\left(f(x),g(x)\right)\equiv0$ on $X$.
$X$ is isomorphic to $X_F$ by extension of
\[X\ni x\mapsto\left(f(x),g(x)\right)\in X_F.\]

Our purpose is to find such irreducible polynomials concretely to given especial compact Riemann surfaces.
Let $\zeta_p$ be a primitive $p$-th root of unity.
The following proposition motivates us to the main theorem of this paper.
\begin{prop}\label{pre1}
Let $X$ be a compact Riemann surface which has an automorphism $\tau$ of order $p$ such that the genus of $X/\langle\tau\rangle$ is zero.
Then an equation of X is given by a semi-hyperelliptic curve
\begin{eqnarray}
y^p=\prod_{i=1}^r(x-a_i)^{m_i}\label{eq:semi}
\end{eqnarray}
and $\tau$ corresponds to $(x,y)\mapsto (x,\zeta_py)$ on this semi-hyperelliptic curve.
\begin{proof}
Let  $\tau^*$ be an automorphism on ${\cal M}(X)$ defined by $\tau^*(f)=f\circ\tau$.
Since
\[{\cal M}(X)^{\langle\tau^*\rangle}=\left\{\,f\in{\cal M}(X)\,|\,\tau^*(f)=f\,\right\}\]
is isomorphic to ${\cal M}(X/\langle\tau\rangle)$, there is a meromorphic function $\mathbf{x}\in{\cal M}(X)$ such that ${\cal M}(X)^{\langle\tau^*\rangle}=\mathbb{C}(\mathbf{x})$.
One can easily see that the degree of $\mathbf{x}$ is $p$ and then $\left[{\cal M}(X):\mathbb{C}(\mathbf{x})\right]=p<\infty$.
We regard the automorphism $\tau^*$ as an endomorphism of the vector space ${\cal M}(X)$ over the field $\mathbb{C}(\mathbf{x})$.
We claim that the minimal polynomial of $\tau^*$ is $t^p-1\in\mathbb{C}(\mathbf{x})[t]$.
Since $(\tau^*)^p={\rm id}$, it suffices to show the minimality of the degree.
We assume that
\begin{eqnarray}
a_0(\mathbf{x})+a_1(\mathbf{x})\,\tau^*+\cdots+a_{p-1}(\mathbf{x})(\tau^*)^{p-1}=0\,,\label{eq:degree}
\end{eqnarray}
where $a_0(\mathbf{x}),\cdots,a_{p-1}(\mathbf{x})\in\mathbb{C}(\mathbf{x})$.
We should show $a_l(\mathbf{x})\equiv0$ for all $l=0,\cdots,p-1$.
We take a point $P\in X$ such that $P,\,\tau(P),\,\cdots,\,\tau^{p-1}(P)$ differ from each other and $a_l(\mathbf{x})(P)$ is finite for all $l$.
We also take meromorphic functions $f_j$ for $j=1,\cdots,p$ such that $f_j\left(\tau^l(P)\right)=j^l$.
The reason of the existence of $f_j$ is shown by the next theorem.
For a proof, see \cite{Fr} for example.

\begin{theo}
Let $X$ be a compact Riemann surface and $S\subset X$ be a finite subset.
Assume that for each $s\in S$ a complex number $a_s\in\mathbb{C}$ is given.
Then there is a meromorphic function $f\in{\cal M}(X)$ such that $f(s)=a_s$ for all $s\in S$.
\end{theo}

By substituting $f_j(P)$ for (\ref{eq:degree}), we have
\[a_0(\mathbf{x})(P)+a_1(\mathbf{x})(P)\cdot j+\cdots+a_{p-1}(\mathbf{x})(P)\cdot j^{p-1}=0.\]
This means that at most $(p-1)$-th degree polynomial has $p$ distinct solution, then we have $a_l(\mathbf{x})(P)=0$.
Since this argument is held for all such points $P\in X$, we have $a_l(\mathbf{x})\equiv0$.
Therefore, the minimal polynomial of $\tau^*$ is $t^p-1$.

By Cayley-Hamilton theorem, $\zeta_p$ is an eigenvalue of $\tau^*$.
We take an eigenvector $\mathbf{y}\in{\cal M}(X)$, that is, $\tau^*(\mathbf{y})=\zeta_p\mathbf{y}$.
By $\tau^*(\mathbf{y}^p)=\mathbf{y}^p$, we get $\displaystyle\mathbf{y}^p=\frac{a(\mathbf{x})}{b(\mathbf{x})}\in \mathbb{C}(\mathbf{x})$.
Replacing $\mathbf{y}$ with $\displaystyle\frac{\mathbf{y}}{b(\mathbf{x})}$, we can assume $\mathbf{y}^p$ is an element of $\mathbb{C}[\mathbf{x}]$.
We can also assume $\mathbf{y}^p$ is a monic polynomial in $\mathbf{x}$ by replacing $\mathbf{y}$ with $c\mathbf{y}$ for a suitable non-zero constant $c$.
Furthermore, since $\mathbf{y}$ is injective on the generic fiber of $\mathbf{x}$, we have ${\cal M}(X)=\mathbb{C}(\mathbf{x},\mathbf{y})$.
Then we conclude that an equation of X is given by the semi-hyperelliptic curve (\ref{eq:semi}).

We finally check the behavior of $\tau$ on this semi-hyperelliptic curve.
Let 
\[\Phi:X\ni P\longmapsto \left(\mathbf{x}(P),\,\mathbf{y}(P)\right)\in\left\{\,{\rm The\ compact\ Riemann\ surface\ given\ by}\ (\ref{eq:semi})\,\right\}\]
be an isomorphism and $(x,y)=(\mathbf{x}(P),\,\mathbf{y}(P))$.
By
\begin{align*}
\Phi\circ\tau\circ\Phi^{-1}(x,y) &= \Phi\circ\tau(P)\\[1mm]
&= \left(\mathbf{x}\circ\tau(P),\,\mathbf{y}\circ\tau(P)\right)\\[1mm]
&= \left(\mathbf{x}(P),\,\zeta_p\mathbf{y}(P)\right)\\[1mm]
&= (x,\zeta_py)
\end{align*}
we have that $\tau$ corresponds to $(x,y)\mapsto (x,\zeta_py)$.
\end{proof}
\end{prop}
\begin{rem}
The irreducibility of
\[y^p-\prod_{i=1}^r(x-a_i)\]
is shown by the behavior of $\tau$ on the semi-hyperelliptic curve.
If it is a reducible polynomial, the map $(x,y)\mapsto(x,\zeta_py)$ is not defined on a compact Riemann surface defined by $y^{p'}\hspace{-3pt}=\prod_{i=1}^{r'}(x-a_i')$, here $p'<p$.
\end{rem}
By repeating the replacement of $y$ to $(x-a_i)y$, we can assume $1\le m_i<p$.
Our next goal is to determine $r,\,m_i$ and $a_i$ in the definition of a semi-hyperelliptic curve.
To do this, we have to consider some properties of a semi-hyperelliptic curve.
Let $X$ be a compact Riemann surface defined by (\ref{eq:semi}) and we see that how $X$ is done its compactification.
By (\ref{eq:com}), $X$ is obtained by the compactification of
\begin{eqnarray}
\left\{\,(x,y)\in\mathbb{C}^2\,|\,y^p=\prod_{i=1}^r(x-a_i)^{m_i},y\ne0\,\right\}\label{eq:add}
\end{eqnarray}
and thus, points which are related to $(a_i,0)$ and infinity points are added.
To see that how to add points to the curve (\ref{eq:add}) by compactification,
we shall define a chart $\varphi_P$ around each point $P$.

Let $\g(a,b)$ be the greatest common divisor of $a$ and $b$.

If $P=(a_i,0)$, we consider
\[\varphi_P^{-1}(t)=
\left(\,t^{\frac{p}{\g(p,m_i)}}+a_i, \,t^{\frac{m_i}{\g(p,m_i)}}\hspace{-3pt}\sqrt[p]{\prod_{k\ne i}\left(t^{\frac{p}{\g(p,m_i)}}+a_i-a_k\right)^{m_k}}\,\right)\]
defined in a small disc.
Then the branch order of the first projection $X\ni(x,y)\mapsto x\in\hat{\mathbb{C}}$ at $(a_i,0)$ is $\frac{p}{\g(p,m_i)}$.
To be the degree of this projection is $p$, we have to add $\g(p,m_i)$ points $(a_{i,1},0),\cdots,(a_{i,\g(p,m_i)},0)$.
We sometimes use simply notation $a_{i,l}$ instead of $(a_{i,l},0)$.
We also use $(a_i,0)$ if $\g(p,m_i)=1$.
We see the branch order of $\tau:(x,\,y)\mapsto(x,\,\zeta_py)$ at $a_{i,l}$ is $\frac{p}{\g(p,m_i)}$ by considering the natural projection $X\rightarrow X/\langle\tau\rangle$.

If $P$ is a infinity point, we also consider
\[\varphi_P^{-1}(t)=
\begin{cases}
\left(\,t^{-\frac{p}{\g(p,m)}}, \,t^{-\frac{m}{\g(p,m)}}\sqrt[p]{\prod_{i=1}^r\left(1-a_it^{\frac{p}{\g(p,m)}}\right)^{m_i}}\ \right)&(\,0<|t|<\varepsilon\,)\\[2mm]
\infty&(\,t=0\,)\hspace{5pt},
\end{cases}\]
where $m$ is $\sum_{i=1}^r m_i$.
Since the branch order of the first projection at a infinity point is $\frac{p}{\g(p,m)}$, 
we also need to add $\g(p,m)$ points $\infty_1,\,\cdots,\,\infty_{\g(p,m)}$ and the branch order of $\tau$ at $\infty_l$ is $\frac{p}{\g(p,m)}$.
In particular, the infinity points are non-branched points of $\tau$ if and only if $m$ is divided by $p$.

In Proposition \ref{pre1}, we take $\Psi$ be an isomorphism from $X/\langle\tau\rangle$ to $\hat{\mathbb{C}}$ and let $Q_1,\cdots,Q_s\in X/\langle\tau\rangle$ be the branched values of the natural projection $X\rightarrow X/\langle\tau\rangle$.
By composing $\Psi$ with a M${\rm \ddot{o}}$bius transformation if necessary,
we can assume $\Psi(Q_i)$ is contained in $\mathbb{C}$ for all $i$.
By the next commutative diagram
\[\begin{CD}
X @>{\rm The\ natural\ projection}>> X/\langle\tau\rangle\\
@V\Phi VV @VV\Psi V\\
X'@>>{\rm \hspace{4pt}The\ first\ projection\hspace{4pt}}> \hat{\mathbb{C}}
\end{CD}\]
where $X'$ is the compactification of (\ref{eq:add}), we have
\[\Psi\left(\left\{Q_1,\cdots,Q_s\right\}\right)=\left\{a_1,\cdots,a_r\right\}.\]
Thus, we can assume $r=s,\Psi(Q_i)=a_i$ and the infinity points are non-branched points of $\tau$.
By the following facts, we obtain conditions about $m_i$.
Under the assumption of Proposition \ref{pre1}, let $P_{1,1},\cdots,P_{1,n_1},\cdots,P_{r,1},\cdots,P_{r,n_r}$ be the branched points of $\tau$,
 where $P_{i,1},\,\cdots,\,P_{i,n_i}$ are $\tau$ equivalent points.
Namely, $\tau(P_{i,1})=P_{i,2},\tau(P_{i,2})=P_{i,3},\cdots,\tau(P_{i,n_i})=P_{i,1}$.
Then the Riemann surface $X$ is given by
\[y^p=\prod_{i=1}^r\left(x-\Psi\left([P_{i,1}]\right)\,\right)^{m_i}\]
and the exponents $m_i$ are satisfied with $\g(p,m_i)=n_i$ for all $i$ and $\sum_{i=1}^rm_i$ is divided by $p$.

In order to completely determine $m_i$, we define the "rotation number" (cf. \cite{GG}).
Let $X$ be a Riemann surface which has an automorphism $\tau$ of order $p$, and a point $P\in X$ be fixed.
Take $n$ to be the smallest natural number such that $\tau^n(P)=P$ and $\varphi$ to be a $\tau$-invariant chart around $P$ centered at the origin.
Then $\varphi\circ\tau^n\circ\varphi^{-1}$ is an automorphism of a small disk fixing the origin with order $\frac{p}{n}$.
Hence, it is of the form
\begin{eqnarray}
\varphi\circ\tau^n\circ\varphi^{-1}(t)=\zeta_{\frac{p}{n}}^k\cdot t\,.\label{eq:rotation}
\end{eqnarray}
Here, $k$ is the integer with $0\le k<\frac{p}{n}$.
We call the pair of $n$ and $k$ a rotation number of $\tau$ at $P$, and we denote it by ${\cal R}_\tau(P)={\cal R}(n,k)$.

\begin{rem}
The number $k$ is independent of the choice of the chart $\varphi$, 
and we see that ${\cal R}_\tau(P)={\cal R}_\tau(P')$ if $P$ and $P'$ are $\tau$ equivalent.
\end{rem}

We often consider the rotation number at branched points of the natural projection $X\rightarrow X/\langle\tau\rangle$.
Actually, we simply have ${\cal R}_\tau(P)={\cal R}(p,0)$ if $P$ is a non-branched point of this projection.
On other hand, if $P$ is a fixed point of $\tau$, we have $n=1$ and the rotation number's concept, for these points, is only the exponent $k$ of (\ref{eq:rotation}).

For getting exponents $m_i$, we consider the rotation number at $a_{i,l}$ of a semi-hyperelliptic curve (\ref{eq:semi}).
\begin{lem}\label{uni}
Let $X$ be a semi-hyperelliptic curve and $\tau$ be an automorphism $(x,\,y)\mapsto (x,\,\zeta_py)$ of $X$ as before.
Let $k$ be a unique number satisfying $k\cdot\frac{m_i}{\g(p,m_i)}\equiv1$ mod $\frac{p}{\g(p,m_i)}$ with $1\le k<\frac{p}{\g(p,m_i)}$.
Then
\[{\cal R}_\tau(a_{i,l})={\cal R}(\,\g(p,m_i),k).\]
\begin{proof}
The existence and uniqueness of $k$ are shown by an elementary argument.
Indeed, if $a$ and $b$ are coprime integers, the equation $ak+bl=1$ has a unique number solution $k$ with $1\le k<b$ for the suitable integer $l$.

It is clear that the smallest number $n$ with $\tau^n(a_{i,l})=a_{i,l}$ is $\g(p,m_i)$ since we add $\g(p,m_i)$ points $a_{i,1},\cdots,{a_{i,\g(p,m_i)}}$ to (\ref{eq:add}).
Then by
\begin{align*}
&\varphi_{a_{i,l}}\circ\tau^{\g(p,m_i)}\circ\varphi_{a_{i,l}}^{-1}(t)\\[2mm]
&=\left(\varphi_{a_{i,l}}\circ\tau^{\g(p,m_i)}\right)\left(\,\left(\,t^{\frac{p}{\g(p,m_i)}}+a_i, \,t^{\frac{m_i}{\g(p,m_i)}}\hspace{-4pt}\sqrt[p]{\prod_{k\ne i}\left(t^{\frac{p}{\g(p,m_i)}}+a_i-a_k\right)^{m_k}}\,\right)\,\right)\\[3mm]
&=
\varphi_{a_{i,l}}\,\left(\,t^{\frac{p}{\g(p,m_i)}}+a_i, \,\zeta_{\frac{p}{\g(p,m_i)}}\hspace{-3pt}\cdot\hspace{-1pt}t^{\frac{m_i}{\g(p,m_i)}}\hspace{-4pt}\sqrt[p]{\prod_{k\ne i}\left(t^{\frac{p}{\g(p,m_i)}}+a_i-a_k\right)^{m_k}}\,\right)\\[3mm]
&=\zeta_{\frac{p}{\g(p,m_i)}}^k\cdot t\ ,
\end{align*}
we have ${\cal R}_\tau(a_{i,l})={\cal R}(\,\g(p,m_i),\,k)$.
\end{proof}
\end{lem}

From above results, we obtain the next theorem.
\begin{theo}\label{main1}
Let $X$ be a compact Riemann surface which has an automorphism $\tau$ of order $p$ such that the genus of $X/\langle\tau\rangle$ is zero.
$P_{1,1},\cdots,P_{1,n_1},\cdots,P_{r,1},\cdots,P_{r,n_r}$ are all branched points of $\tau$, where $P_{i,1},\cdots,P_{i,n_i}$ are $\tau$ equivalent points and the rotation numbers are given by ${\cal R}_\tau(a_{i,l})={\cal R}(n_i,k_i)$ for all $i$.
If we take unique numbers $m_i$ such that $n_i=\g(p,m_i)$ and $k_i\cdot\frac{m_i}{n_i}\equiv1$ mod $\frac{p}{n_i}$ with $1\le m_i<p$, 
an equation of X is a semi-hyperelliptic curve given by
\[y^p=\prod_{i=1}^r(x-a_i)^{m_i}\]
and $\tau$ corresponds to $(x,y)\mapsto (x,\zeta_py)$.
Furthermore, if we take an isomorphism $\Psi$ from $X/\langle\tau\rangle$ to $\hat{\mathbb{C}}$, $a_i$ is given by $\Psi\left([P_{i,1}]\right)$.
\end{theo}

The existence and uniqueness of $m_i$ are shown by the same argument as in the proof of Lemma \ref{uni}.
Of course, $\sum_{i=1}^rm_i$ is divided by $p$ in this case.
Finally, we give a remark about normalization. 
By composing a suitable M${\rm \ddot{o}}$bius transformation to $\Psi$, we get a normalized equation of $X$, namely
\[y^p=x^{m_1}(x-1)^{m_2}(x-a_3')^{m_3}\cdots(x-a_{r-1}')^{m_{r-1}}.\]
Therefore, if $r\le3$, we immediately determine an equation completely, which corresponds to $X_q$ for $q\le7$.
\section{The genus of quotient compact Riemann surfaces of $X_q$}
In this section, we see that there is an automorphism $\tau$ of $X_q$ such that the genus of $X_q/\langle\tau\rangle$ is zero for $q\le10$ or $q=12$.\vspace{2pt}

Let $\hat{\mathbb{H}}$ be $\mathbb{H}\cup\mathbb{Q}\cup\{\infty\}$, and we give a unique topology of $\hat{\mathbb{H}}$ such that it satisfies the following properties.\\
\hspace{10pt}1)\hspace{5pt}
The topology induced by $\hat{\mathbb{H}}$ gives the usual topology on $\mathbb{H}$.\\
\hspace{10pt}2)\hspace{5pt}
Elements of $SL\left(2,\mathbb{Z}\right)$ acts continuously on $\hat{\mathbb{H}}$.\\
\hspace{10pt}3)\hspace{5pt}
A subset of $\hat{\mathbb{H}}$ is a neighborhood of $\infty$ if and only if it contains a set $\left\{\,z\in\mathbb{H}:{\rm Im}\,z>C\,\right\}\cup\{\infty\}$\\
\hspace{22pt}for a positive number $C>0$.\\[-2.5mm]

Since $X_q=\overline{\mathbb{H}/\Gamma_q}$ is isomorphic to $\hat{\mathbb{H}}/\Gamma_q$, we redefined $X_q$ by $\hat{\mathbb{H}}/\Gamma_q$.
We take an automorphism $\tau_n$ of $X_q$ given by $[z]\mapsto[z+n]$, where $n$ is a positive divisor of $q$.
It is grad if the genus of $X_q/\langle\tau_n\rangle$ is zero.
\begin{rem}
We naturally think that the genus of $X_q/\langle\tau\rangle$ decreases as the order of $\tau\in{\rm Aut}(X_q)$ increases \footnote{Of coures, there are opposite cases. For example, see table 7 in appendix.}.
The order of automorphisms of $X_q$ is at most $q$ for $q\ge7$, and the remainder of $q$ divided by 4 isn't 2 or $q$ is divided by 3 (see Proposition \ref{most} in appendix).
For example, $q=7,8,9,11,12,13,15,\cdots$.
The order of the automorphism $\tau_1:[z]\mapsto[z+1]$ reaches the bound for these $q$.
Thus, taking $\tau_n$, especially $\tau_1$, from automorphisms of $X_q$ is reasonable.
\end{rem}
For getting the genus of $X_q/\langle\tau_n\rangle$, we define some notations first.
\begin{nota}
For $\Gamma\subset SL\left(2,\mathbb{Z}\right)$, we define $\tilde{\Gamma}$ by $\Gamma\cup-\Gamma$ and for $\hat{\mathbb{Q}}=\mathbb{Q}\cup\{\infty\}$, we define next notations.
\[S_q\hspace{-2pt}:=\hat{\mathbb{Q}}/\Gamma_q\,,\hspace{3pt}h_q\hspace{-2pt}:=\#S_q\,,\hspace{3pt}R_q\hspace{-2pt}:=\left[SL\left(2,\mathbb{Z}\right):\tilde{\Gamma}_q\right]\]
Let $X_q^n$ be $\hat{\mathbb{H}}/\Gamma_q^n$, which is isomorphic to $X_q/\langle\tau_n\rangle$, where
\[\Gamma_q^n:=\left\{\begin{pmatrix}
a&b\\c&d
\end{pmatrix}
\in SL\left(2,\mathbb{Z}\right)\,|\,
a\equiv d\equiv1,\,c\equiv0\,\,({\rm mod}\ q),\,b\equiv0\,\,({\rm mod}\ n)\right\}.\]
We also define
\[S_q^n\hspace{-2pt}:=\hat{\mathbb{Q}}/\Gamma_q^n\,,\hspace{3pt}h_q^n\hspace{-2pt}:=\#S_q^n\,,\hspace{3pt}R_q^n\hspace{-2pt}:=\left[SL\left(2,\mathbb{Z}\right):\tilde{\Gamma}_q^n\right]\,.\]
We set $g_q$ and $g_q^n$ be the genera of $X_q$ and $X_q^n$, respectively.

Finally, for a natural number $q\in\mathbb{N}$, let ${\cal P}(q)$ be the set consisting of all primes $l$ which divides $q$.
For example, ${\cal P}(12)=\{2,3\}$.
\end{nota}

The next theorem is well known.
For a proof, see \cite{Fr} or \cite{Si} for example.
\begin{theo}\label{fre}
For $q\ge3$, we have
\[R_q=\frac{q^3}{2}\pp\hspace{2pt},\hspace{3pt}
h_q=\frac{q^2}{2}\pp\hspace{2pt},\hspace{3pt}
g_q=1+\frac{(q-6)q^2}{24}\pp\hspace{3pt}.\]
\end{theo}

We remark that $R_q=\#PSL\left(2,\mathbb{Z}/q\mathbb{Z}\right)$.
Before getting $g_q^n$, we evaluate $R_q^n$ and $h_q^n$.
It is because we have $g_q^n=1-\frac{h_q^n}{2}+\frac{R_q^n}{12}$ if $\Gamma_q^n$ acts freely on $\mathbb{H}$.
Actually, $\Gamma_q^n$ acts freely on $\mathbb{H}$ for $q\ge4$.
\begin{prop}
For $q\ge3$, we have
\[R_q^n=\frac{nq^2}{2}\pp.\]
\begin{proof}
We get it by \ $\left[SL\left(2,\mathbb{Z}\right):\tilde{\Gamma}_q\right]=\left[SL\left(2,\mathbb{Z}\right):\tilde{\Gamma}_q^n\right]\cdot\left[\tilde{\Gamma}_q^n:\tilde{\Gamma}_q\right]$ 
\ and \ $\left[\tilde{\Gamma}_q^n:\tilde{\Gamma}_q\right]=\left[\Gamma_q^n:\Gamma_q\right]=\frac{q}{n}$.
\end{proof}
\end{prop}

Then we consider $h_q^n$.
Let $\Gamma$ be a subgroup of finite index of $SL\left(2,\mathbb{Z}\right)$ and $\kappa$ be an element of $\hat{\mathbb{Q}}$. 
We take $N\in SL\left(2,\mathbb{Z}\right)$ such that $N(\infty)=\kappa$.
Then there is a positive number $R$ such that
\[\left\{\,M\in N^{-1}\tilde{\Gamma} N\,|\,M(\infty)=\infty\,\right\}
=\left\{\,\pm\begin{pmatrix}
1&mR\\0&1
\end{pmatrix}
|\,m\in\mathbb{Z}\,\right\}\,.\]
We call $R$ the width of $\kappa$ and use the notation ${\cal W}_\Gamma(\kappa)$ \cite{Fr}.
\begin{rem}
This definition is independent of the choice of $N$ since if we take another $N'$, it satisfies 
$N'=N\left(\begin{smallmatrix}
1&b\\0&1
\end{smallmatrix}\right)$ for some $b\in\mathbb{Z}$.
Moreover, it depends only on the $\Gamma$-equivalence class.
It is because if $\kappa$ and $\kappa'$ are satisfied $\gamma(\kappa)=\kappa'$ for some $\gamma\in\Gamma$, we have $\gamma N(\infty)=\kappa'$ and then $N^{-1}\tilde{\Gamma} N=(\gamma N)^{-1}\tilde{\Gamma}(\gamma N)$.
Therefore, we can define the width of elements of $\hat{\mathbb{Q}}/\Gamma$ in a natural way and we use the same notation.
\end{rem}
\begin{lem}\label{width}
Let $\Gamma$ and $\Gamma'$ be subgroups of finite index of $SL\left(2,\mathbb{Z}\right)$ and each of them contains the negative unit matrix.
Set $\Gamma$ to be a subgroup of $\Gamma'$ and $p$ to be $\left[\Gamma':\Gamma\right]$\,.
Let us take an element from $\hat{\mathbb{Q}}/\Gamma'$ and let $\kappa$ denote its representative.
We consider the quotient of the stabilizer of $\kappa$ in $\Gamma'$ determined by $\Gamma$.
Let $\kappa_1,\cdots,\kappa_h$ denote representatives of this quotient's elements.
Then
\[p\cdot{\cal W}_{\Gamma'}(\kappa)=\sum_{i=1}^h{\cal W}_\Gamma(\kappa_i)\,.\]
\begin{proof}
Let $SL\left(2,\mathbb{Z}\right)_\kappa$ and $\Gamma_\kappa'$ be subgroups of $SL\left(2,\mathbb{Z}\right)$ and $\Gamma'$, respectively.
Each of them fixes $\kappa$.
Since ${\cal W}_{\Gamma'}(\kappa)=\left[SL\left(2,\mathbb{Z}\right)_\kappa:\Gamma_\kappa'\right]$, we take $N_1,\cdots,N_{{\cal W}_{\Gamma'}(\kappa)}$ which are the set of left cosets of $\Gamma'_\kappa$ in $SL\left(2,\mathbb{Z}\right)_\kappa$.
For $i=1,\cdots,h$, let $SL\left(2,\mathbb{Z}\right)_{\kappa_i}$ and $\Gamma_{\kappa_i}$ be subgroups of $SL\left(2,\mathbb{Z}\right)$ and $\Gamma$ such that they fix $\kappa_i$, respectively.
Since ${\cal W}_{\Gamma}(\kappa_i)=\left[SL\left(2,\mathbb{Z}\right)_{\kappa_i}:\Gamma_{\kappa_i}\right]$, we also take $N_{i,1},\cdots,N_{i,{\cal W}_{\Gamma}(\kappa_i)}$ which are the set of left cosets of $\Gamma_{\kappa_i}$ in $SL\left(2,\mathbb{Z}\right)_{\kappa_i}$.
Then by letting $A_1,\cdots,A_p$ denote the set of left cosets of $\Gamma$ in $\Gamma'$ and 
$M,M_1,\cdots,M_h$ satisfy $M(\kappa)=M_1(\kappa_1)=\cdots=M_h(\kappa_h)=\infty$, 
we claim that
\begin{eqnarray}
\bigoplus_{\mu,\,\nu}MN_\mu A_\nu\Gamma=\bigoplus_{i,j}M_iN_{i,j}\Gamma\,.\label{eq:oplus}
\end{eqnarray}
We prove (\ref{eq:oplus}) in several steps.\\[1mm]
{\it Step 1}: $MN_\mu A_\nu\Gamma\cap MN_{\mu'}A_{\nu'}\Gamma=\phi$.\\
We assume $MN_\mu A_\nu\Gamma\cap MN_{\mu'}A_{\nu'}\Gamma\ne\phi$.
We should show $\mu=\mu'$ and $\nu=\nu'$.
There is $\gamma\in\Gamma$ such that
\[MN_\mu A_\nu=MN_{\mu'}A_{\nu'}\gamma\Leftrightarrow N^{-1}_{\mu'}N_\mu=A_{\nu'}\gamma A^{-1}_\nu\,.\]
Since $A_{\nu'}\gamma A^{-1}_\nu\in\Gamma'$, we get $N^{-1}_{\mu'}N_\mu\in\Gamma'_\kappa$ and so $\mu=\mu'$.
We also get $\nu=\nu'$ by $A_\nu\Gamma\cap A_{\nu'}\Gamma\ne\phi$.\\[1mm]
{\it Step 2}: $M_iN_{i,j}\Gamma\cap M_{i'}N_{i',j'}\Gamma=\phi$.\\
We assume $M_iN_{i,j}\Gamma\cap M_{i'}N_{i',j'}\Gamma\ne\phi$.
We should show $i=i'$ and $j=j'$.
There is $\gamma\in\Gamma$ such that
\[M_iN_{i,j}=M_{i'}N_{i',j'}\gamma\Leftrightarrow N^{-1}_{i',j'}M^{-1}_{i'}M_iN_{i,j}=\gamma\,.\]
Then we get
\[\gamma(\kappa_i)=N^{-1}_{i',j'}M^{-1}_{i'}M_iN_{i,j}(\kappa_i)=N^{-1}_{i',j'}M^{-1}_{i'}(\infty)=\kappa_{i'}\]
and so $i=i'$.
Since $N_{i,j}\Gamma\cap N_{i,j'}\Gamma\ne\phi$, we also get $N_{i,j}\Gamma_{\kappa_i}\cap N_{i,j'}\Gamma_{\kappa_i}\ne\phi$ and so $j=j'$.\\[1mm]
{\it Step 3: The left-hand side is contained in the right-hand side.}\\
It is sufficient to prove that $MN_\mu A_\nu$ is contained in the right-hand side.
We take $\gamma\in\Gamma$ and $i$ such that $A^{-1}_\nu(\kappa)=\gamma(\kappa_i)$.
Since $M^{-1}_iMN_\mu A_\nu\gamma\in SL(2,\mathbb{Z})_{\kappa_i}$, there is $j$ such that $M^{-1}_iMN_\mu A_\nu\gamma\in N_{i,j}\Gamma_{\kappa_i}\subset N_{i,j}\Gamma$.
Then we have $MN_\mu A_\nu\in M_iN_{i,j}\Gamma$.\\[1mm]
{\it Step 4: The right-hand side is contained in the left-hand side.}\\
It is sufficient to prove that $M_iN_{i,j}$ is contained in the left-hand side.
We take $\gamma\in\Gamma$ such that $\gamma(\kappa)=\kappa_i$.
Since $M^{-1}M_iN_{i,j}\gamma\in SL(2,\mathbb{Z})_\kappa$, there is $\mu$ such that $M^{-1}M_iN_{i,j}\gamma\in N_\mu\Gamma'_\kappa\subset N_\mu\Gamma'$.
Then there is $\nu$ such that $M_iN_{i,j}\in MN_\mu A_\nu\Gamma$.

Thus, the equation (\ref{eq:oplus}) is shown and the proof is completed.
\end{proof}
\end{lem}

In special case $\Gamma'=SL\left(2,\mathbb{Z}\right)$, we get
\begin{cor}
Let $\Gamma$ be a subgroup of finite index of $SL\left(2,\mathbb{Z}\right)$.
Then
\[\left[SL\left(2,\mathbb{Z}\right):\tilde{\Gamma}\right]=\sum_{\kappa\in\hat{\mathbb{Q}}/\Gamma}{\cal W}_\Gamma(\kappa)\,.\]
\end{cor}

We then evaluate $h_q^n$ by using a width.
\begin{nota}
Let $p$ be a natural number and $\prod_{i=1}^kp_i^{r_i}$ be the prime factorization of it.
We define a multiplicative function ${\cal N}(p)$ by
\[{\cal N}(p):=\prod_{i=1}^k\left(1+r_i\cdot\frac{p_i-1}{p_i+1}\right)\,.\]
Here, we note that ${\cal N}(1)$ is 1.
\end{nota}
\begin{prop}\label{h}
For $q\ge5$, we have
\[h_q^n=\frac{nq\cdot{\cal N}\hspace{-1pt}\left(\frac{q}{n}\right)}{2}\pp\,.\]
\begin{proof}
Let $p$ be $\frac{q}{n}$.
We split the proof into several steps.\\
{\it Step 1}:
We describe by $\frac{x}{z}$ a representative of $\kappa\in S_q^n$, where $x$ and $z$ are coprime integers.
We claim that
\[{\cal W}_{\Gamma_q^n}(\kappa)=\frac{q}{\g(p,z)}\,.\]
Here, we note that $\infty$ is $\frac{1}{0}$ and $\g(p,0)$ is defined to be $p$.

Let $k$ be $\g(p,z)$ and $p',z'$ are coprime integers such that $p=kp',z=kz'$.
We take 
$N=\left(\begin{smallmatrix}
x&y\\z&w
\end{smallmatrix}\right)\in SL\left(2,\mathbb{Z}\right)$ 
and consider elements of $N^{-1}\Gamma_q^nN$ such that $\infty$ is fixed.
By
\[\begin{pmatrix}w&-y\\-z&x\end{pmatrix}
\begin{pmatrix}-qxz'R+1&np'x^2R\\-qzz'R&qxz'R+1\end{pmatrix}
\begin{pmatrix}x&y\\z&w\end{pmatrix}
=\begin{pmatrix}
1&np'R\\0&1
\end{pmatrix}\]
we have ${\cal W}_{\Gamma_q^n}(\kappa)\le np'$.
Then we should show ${\cal W}_{\Gamma_q^n}(\kappa)\ge np'$.

We consider all the elements in $N^{-1}\Gamma_q^nN$ which fix $\infty$.
Since
\begin{eqnarray}\label{eq:step1}
\begin{pmatrix}w&-y\\-z&x\end{pmatrix}
\begin{pmatrix}1&nb\\0&1\end{pmatrix}
\begin{pmatrix}x&y\\z&w\end{pmatrix}
\equiv\begin{pmatrix}
nzwb+1&nw^2b\\-nz^2b&-nzwb+1
\end{pmatrix}\hspace{10pt}({\rm mod}\ q),
\end{eqnarray}
we see that $nzwb$ is a multiple of $q$, that is, $wb$ is divided by $p'$.
If not, $nzwb\equiv-2\equiv2$ in modular $q$ by diagonal components of (\ref{eq:step1}).
It is a contradiction with $q\ge5$.
Thus, by (1,2) component of (\ref{eq:step1}), we have ${\cal W}_{\Gamma_q^n}(\kappa)\ge np'$.
\begin{rem}
The condition $q\ge5$ in Proposition \ref{h} is owing to Step 1.
Indeed, we have
\begin{align*}
\begin{pmatrix}1&0\\2&1\end{pmatrix}^{-1}
\begin{pmatrix}-4m+1&2m\\-8m&4m+1\end{pmatrix}
\begin{pmatrix}1&0\\2&1\end{pmatrix}
&=\begin{pmatrix}1&2m\\0&1\end{pmatrix}\\[2mm]
\begin{pmatrix}1&0\\2&1\end{pmatrix}^{-1}
\begin{pmatrix}-4m-1&2m+1\\-8m-4&4m+3\end{pmatrix}
\begin{pmatrix}1&0\\2&1\end{pmatrix}
&=\begin{pmatrix}1&2m+1\\0&1\end{pmatrix}.
\end{align*}
Thus, we obtain
\[\left\{M\in\begin{pmatrix}1&0\\2&1\end{pmatrix}^{-1}
\tilde{\Gamma}_4^1\begin{pmatrix}1&0\\2&1\end{pmatrix}|\,
M(\infty)=\infty\right\}
=\left\{\pm\begin{pmatrix}
1&m\\0&1
\end{pmatrix}
|m\in\mathbb{Z}\right\}.\]
This implies ${\cal W}_{\Gamma_4^1}\hspace{-1.5pt}\left(\frac{1}{2}\right)=1\ne2=\frac{4}{\g(4,2)}$.
\end{rem}
\hspace{-15pt}{\it Step 2}:
Let $\prod_{i=1}^kp_i^{r_i}$ be the prime factorization of $p$.
For $p$ and $0\le j_i\le r_i$, we define
\[{\cal N}_1(j_i)={\cal N}_1(p,p_i^{j_i}):=
\begin{cases}
1&(\,j_i=0\,)\\
(p_i+1)p_i^{j_i-1}&(\,1\le j_i\le r_i\,)\,.
\end{cases}\]
We claim that the number of elements of $S_q$ such that the denominators of their representative is divided by $\prod_{i=1}^kp_i^{j_i}$ is
\[\frac{h_q}{\prod_{i=1}^k{\cal N}_1(j_i)}\,.\]

Let $H_m$ be the subgroup of $SL\left(2,\mathbb{Z}\right)$ such that these (2,1) entries are dividid by natural number $m$.
It is sufficient to prove that
\[\left[SL\left(2,\mathbb{Z}\right):H_{\prod_{i=1}^kp_i^{j_i}}\right]=\prod_{i=1}^k{\cal N}_1(j_i)\,.\]
We get it by
\[SL\left(2,\mathbb{Z}\right)=H_{p_1}\oplus\ 
\left(\bigoplus_{l=0}^{p_1-1}\begin{pmatrix}l&-1\\1&0\end{pmatrix}H_{p_1}\right)
\hspace{2pt},\hspace{5pt}
H_{p_1}=\bigoplus_{l=0}^{p_1^{j_1-1}-1}
\begin{pmatrix}1&0\\lp_1&1\end{pmatrix}H_{p_1^{j_1}}\]
and so on.\\[1mm] 
{\it Step 3}:
For $p$, we define
\[{\cal N}_2(j_i)={\cal N}_2(p,p_i^{j_i}):=
\begin{cases}
\frac{p_i}{p_i+1}&(\,j_i=0\,)\\[2mm]
\frac{(p_i-1)p_i^{j_i}}{p_i+1}&(\,1\le j_i\le r_i-1\,)\\[2mm]
\frac{p_i^{r_i+1}}{p_i+1}&(\,j_i=r_i\,)\,.
\end{cases}\]
Let $A$, whose width of its element with respect to $\Gamma_q^n$ \footnote{By taking a representative of elements of $S_q$, we may consider the its width with respect to $\Gamma_q^n$ by a natural way.} is $n\prod_{i=1}^kp_i^{j_i}$, be the subset of $S_q$.
We claim that
\[\# A=\frac{h_q}{p}\prod_{i=1}^k{\cal N}_2(j_i)\,.\]

For elements of $A$, by Step 1, the denominators of representatives are divided by $\prod_{i=1}^kp_i^{r_i-j_i}$ and are not divided by $p_s^{r_s-j_s+1}\prod_{i\ne s}p_i^{r_i-j_i}$ for all $s=s_1,\,\cdots,\,s_\mu$ and for positive $j_s$.
By Step 2, we get
\begin{align}\label{eq:many}
\# A&=
\frac{h_N}{\prod_{i=1}^k{\cal N}_1(r_i-j_i)}
-\sum_{s_\nu}\frac{h_N}{{\cal N}_1(r_s-j_s+1)\prod_{i\ne s}{\cal N}_1(r_i-j_i)}\notag \\[1mm]
&\hspace{25pt}+\sum_{s_{\nu},\,s_{\nu'}}\frac{h_N}{{\cal N}_1(r_{s_1}-j_{s_1}+1){\cal N}_1(r_{s_2}-j_{s_2}+1)\prod_{i\ne s_1,s_2}{\cal N}_1(r_i-j_i)}\notag \\[1mm]
&\hspace{25pt}-\cdots+(-1)^\mu\sum_{s_1,\dots,s_\mu}\frac{h_N}{\prod_{\nu=1}^\mu{\cal N}_1(r_{s_\nu}-j_{s_\nu}+1)\,\cdot\prod_{i\ne s_1,\cdots,s_\mu }{\cal N}_1(r_i-j_i)}\notag \\
&=\frac{h_N}{p}\prod_{i=1}^k{\cal N}_2(j_i)\,.
\end{align}
The last equality is showed as follows.
We assume that $n\prod_{i=1}^kp_i^{j_i}=np_{h+1}^{j_{h+1}}\cdots p_{h+l}^{j_{h+l}}\cdot p_{h+l+1}^{r_{h+l+1}}\cdots p_k^{r_k}$ and 
$1\le j_{h+t}<r_{h+t}$ for all $t=1,\cdots,l$.
{\small
\begin{align*}
&\#A\\
&=\frac{h_q}{(p_1+1)p_1^{r_1-1}\cdots(p_h+1)p_h^{r_h-1}(p_{h+1}+1)p_{h+1}^{r_{h+1}-j_{h+1}-1}\cdots(p_{h+l}+1)p_{h+l}^{r_{h+l}-j_{h+l}-1}}\\[2mm]
&\hspace{25pt}-\frac{h_q}{(p_1\p1)p_1^{r_1-1}\cdots(p_h\p1)p_h^{r_h-1}(p_{h+1}\p1)p_{h+1}^{r_{h+1}-i_{h+1}}(p_{h+2}\p1)p_{h+2}^{r_{h+2}-j_{h+2}-1}\cdots(p_{h+l}\p1)p_{h+l}^{r_{h+l}-j_{h+l}-1}}\\[2mm]
&\hspace{25pt}-\cdots
-\frac{h_q}{(p_1\p1)p_1^{r_1\m1}\hspace{-8pt}\cdots(p_h\p1)p_h^{r_h\m1}(p_{h+1}\p1)p_{h+1}^{r_{h+1}\m j_{h+1}\m1}\hspace{-8pt}\cdots(p_{h+l\m1}\p1)p_{h+l\m1}^{r_{h+l\m1}-j_{h+l\m1}\m1}(p_{h+l}\p1)p_{h+l}^{r_{h+l}\m j_{h+l}}}\\[2mm]
&\hspace{25pt}
-\frac{h_q}{(p_1+1)p_1^{r_1-1}\cdots(p_h+1)p_h^{r_h-1}(p_{h+1}+1)p_{h+1}^{r_{h+1}-j_{h+1}-1}\cdots(p_{h+l}+1)p_{h+l}^{r_{h+l}-j_{h+l}-1}(p_{h+l+1}+1)}\\[2mm]
&\hspace{25pt}-\cdots
-\frac{h_q}{(p_1+1)p_1^{r_1-1}\cdots(p_h+1)p_h^{r_h-1}(p_{h+1}+1)p_{h+1}^{r_{h+1}-j_{h+1}-1}\cdots(p_{h+l}+1)p_{h+l}^{r_{h+l}-j_{h+l}-1}(p_k+1)}\\[1.5mm]
&\hspace{25pt}+\cdots\\[1.5mm]
&\hspace{25pt}+\,(-1)^{k-h}
\hspace{-2pt}\frac{h_q}{(p_1\p1)p_1^{r_1-1}\hspace{-8pt}\cdots(p_h\p1)p_h^{r_h-1}(p_{h+1}\p1)p_{h+1}^{r_{h+1}-j_{h+1}}\hspace{-8pt}\cdots(p_{h+l}\p1)p_{h+l}^{r_{h+l}-j_{h+l}}(p_{h+l+1}\p1)\cdots(p_k\p1)}\\[5mm]
&=
\frac{h_q}{p}\cdot\frac{p_1\cdots p_h\,p_{h+1}^{j_{h+1}+1}\cdots p_{h+l}^{j_{h+l}+1}\,p_{h+l+1}^{r_{h+l+1}}\cdots p_k^{r_k}}{(p_1+1)\cdots(p_{h+l}+1)}
-\frac{h_q}{p}\cdot\frac{p_1\cdots p_h\,p_{h+1}^{j_{h+1}}p_{h+2}^{j_{h+2}+1}\cdots p_{h+l}^{j_{h+l}+1}\,p_{h+l+1}^{r_{h+l+1}}\cdots p_k^{r_k}}{(p_1+1)\cdots(p_{h+l}+1)}\\[2mm]
&\hspace{25pt}
-\cdots-\frac{h_q}{p}\cdot\frac{p_1\cdots p_h\,p_{h+1}^{j_{h+1}+1}\cdots p_{h+l-1}^{j_{h+l-1}+1}p_{h+l}^{j_{h+l}}\,p_{h+l+1}^{r_{h+l+1}}\cdots p_k^{r_k}}{(p_1+1)\cdots(p_{h+l}+1)}\\[2mm]
&\hspace{25pt}
-\frac{h_q}{p}\cdot\frac{p_1\cdots p_h\,p_{h+1}^{j_{h+1}+1}\hspace{-5pt}\cdots p_{h+l}^{j_{h+l}+1}\,p_{h+l+1}^{r_{h+l+1}}\cdots p_k^{r_k}}{(p_1+1)\cdots(p_{h+l}+1)(p_{h+l+1}+1)}
-\cdots-\frac{h_q}{p}\cdot\frac{p_1\cdots p_h\,p_{h+1}^{j_{h+1}+1}\hspace{-5pt}\cdots p_{h+l}^{j_{h+l}+1}\,p_{h+l+1}^{r_{h+l+1}}\cdots p_k^{r_k}}{(p_1+1)\cdots(p_{h+l}+1)(p_k+1)}\\[1.5mm]
&\hspace{25pt}+\cdots\\[1mm]
&\hspace{25pt}+\,(-1)^{k-h}
\frac{h_q}{p}\cdot\frac{p_1\cdots p_h\,p_{h+1}^{j_{h+1}}\cdots p_{h+l}^{j_{h+l}}\,p_{h+l+1}^{r_{h+l+1}}\cdots p_k^{r_k}}{(p_1+1)\cdots(p_k+1)}\\[4mm]
&=
\frac{h_q}{p}\cdot\frac{p_1\cdots p_h\,p_{h+1}^{j_{h+1}}\cdots p_{h+l}^{j_{h+l}}\,p_{h+l+1}^{r_{h+l+1}}\cdots p_k^{r_k}}{(p_1+1)\cdots(p_k+1)}\cdot
\big\{\,p_{h+1}\cdots p_{h+l}\,(p_{h+l+1}+1)\cdots(p_k+1)\\[1mm]
&\hspace{25pt}-p_{h+2}\cdots p_{h+l}\,(p_{h+l+1}+1)\cdots(p_k+1)
-\cdots-p_{h+1}\cdots p_{h+l-1}\,(p_{h+l+1}+1)\cdots(p_k+1)\\[1mm]
&\hspace{25pt}-p_{h+1}\cdots p_{h+l}\,(p_{h+l+2}+1)\cdots(p_k+1)
-\cdots-p_{h+1}\cdots p_{h+l}\,(p_{h+l+1}+1)\cdots(p_{k-1}+1)\\[1mm]
&\hspace{25pt}+\cdots+(-1)^{k-h}\,\big\}\\[1mm]
&=
\frac{h_q}{p}\cdot\frac{p_1\cdots p_h\,p_{h+1}^{j_{h+1}}\cdots p_{h+l}^{j_{h+l}}\,p_{h+l+1}^{r_{h+l+1}}\cdots p_k^{r_k}}{(p_1+1)\cdots(p_k+1)}\cdot
(p_{h+1}-1)\cdots(p_{h+l}-1)\cdot(p_{h+l+1}+1-1)\cdots(p_k+1-1)\\
&=
\frac{h_q}{p}\prod_{i=1}^k{\cal N}_2(j_i)
\end{align*}
}

\begin{ex}
Let $q$ be a semiprime number $p_1p_2$.
We consider elements of $S_{p_1p_2}$ such that their width is $p_1p_2$ with respect to $\Gamma_{p_1p_2}^1$.
By Step 1, the width of $\kappa=\left[\frac{x}{z}\right]\in S_{p_1p_2}$ is $p_1p_2$ with respect to $\Gamma_{p_1p_2}^1$ if and only if $\g(p_1p_2,z)=1$.
The denominator $z$ is not divided $p_1$ and $p_2$.
By Step 2 and the following calculation, we count the number of such $\kappa$:
\begin{align*}
&h_{p_1p_2}-\frac{h_{p_1p_2}}{{\cal N}_1(p_1p_2,p_1)}-\frac{h_{p_1p_2}}{{\cal N}_1(p_1p_2,p_2)}+\frac{h_{p_1p_2}}{{\cal N}_1(p_1p_2,p_1)\cdot{\cal N}_1(p_1p_2,p_2)}\\[1mm]
=&h_{p_1p_2}-\frac{h_{p_1p_2}}{p_1+1}-\frac{h_{p_1p_2}}{p_2+1}+\frac{h_{p_1p_2}}{(p_1+1)(p_2+1)}\\[1mm]
=&\frac{h_{p_1p_2}\cdot p_1p_2}{(p_1+1)(p_2+1)}\\[1mm]
=&\frac{h_{p_1p_2}}{p_1p_2}{\cal N}_2(p_1p_2,p_1)\cdot{\cal N}_2(p_1p_2,p_2).
\end{align*}
This equation is corresponding to (\ref{eq:many}) in Step 3.
\end{ex}
\hspace{-15pt}{\it Step 4}:
For $p$, we define
\[{\cal N}_3(j_i)={\cal N}_3(p,p_i^{j_i}):=
\begin{cases}
\frac{p_i}{p_i+1}&(\,j_i=0,\,r_i\,)\\[2mm]
\frac{p_i-1}{p_i+1}&(\,1\le j_i\le r_i-1\,)\,.
\end{cases}\]
Let $B$, whose width of its element is $n\prod_{i=1}^kp_i^{j_i}$, be the subset of $S_q^n$.
We claim that
\[\# B=\frac{h_q}{p}\prod_{i=1}^k{\cal N}_3(j_i)\,.\]

It is shown by Lemma \ref{width}.
Indeed, we have
\[p\cdot n\prod_{i=1}^kp_i^{j_i}\cdot\# B=q\cdot\# A
\ \Leftrightarrow\ 
\# B=\frac{h_q}{p}\prod_{i=1}^k{\cal N}_3(j_i)\,.\]\\[-2mm]
{\it Step 5}:
In this final step, we claim that
\[h_q^n=\sum_{j_1,\cdots,j_k}\frac{h_q}{p}\prod_{i=1}^k{\cal N}_3(j_i)=\frac{nq\cdot{\cal N}(p)}{2}\pp\]
and complete the proof.

It is sufficient to prove that
\[\sum_{j_1,\cdots,j_k}\prod_{i=1}^k{\cal N}_3(j_i)={\cal N}(p)\]
and we show it by induction on $k$.
For $k=1$, we get
\begin{align*}
\sum_{j=0}^r{\cal N}_3(j)&=2\cdot\frac{p}{p+1}+(r-1)\cdot\frac{p-1}{p+1}\\
&=1+r\cdot\frac{p-1}{p+1}\\
&={\cal N}(p)\hspace{5pt}.\\
\intertext{We assume that the claim is held for $k-1$.
Then we have}
\sum_{j_1,\cdots,j_k}\prod_{i=1}^k{\cal N}_3(j_i)&=
\sum_{j_k=0}^{r_k}\,\left(\,\sum_{j_1,\cdots,j_{k-1}}\prod_{i=1}^{k-1}{\cal N}_3(j_i)\right)\cdot{\cal N}_3(j_k)\\[2mm]
&=\sum_{j_k=0}^{r_k}\,{\cal N}\left(\prod_{i=1}^{k-1}p_i^{r_i}\right)\cdot{\cal N}_3(j_k)\\[2mm]
&={\cal N}\left(\prod_{i=1}^{k-1}p_i^{r_i}\right)\cdot{\cal N}(p_k^{r_k})\\[2mm]
&={\cal N}(p)
\end{align*}
and the proof.
\end{proof}
\end{prop}

We recall that $g_q^n=1-\frac{h_q^n}{2}+\frac{R_q^n}{12}$ and then we obtain the next theorem and table 1.
\begin{theo}\label{p}
For $q\ge5$, we have
\[g_q^n=1+\frac{\left(q-6{\cal N}\hspace{-1pt}\left(\frac{q}{n}\right)\right)nq}{24}\pp\,.\]
\end{theo}
\begin{table}[htb]
\begin{center}
\begin{tabular}{c||c|c|c|c|c|c|c|c|c|c|c|c|c|c|c|c|c}
$q$&$1\sim5$&6&7&8&9&10&11&12&13&14&15&16&17&18&19&20&$\cdots$\\[1mm] \hline \hline
$g_q$&0&1&3&5&10&13&26&25&50&49&73&81&133&109&196&169&$\cdots$\\[1mm] \hline
$g_q^1$&0&0&0&0&0&0&1&0&2&1&1&2&5&2&7&3&$\cdots$
\end{tabular}
\caption{Genera of $X_q$ and $X_q^1$}
\end{center}
\end{table}
By considering whether $q-6{\cal N}(q)$ is negative or by table 1, we have the next corollary .
\begin{cor}
For $q\le10$ or $q=12$, $X_q$ has an automorphisms $\tau$ such that the genus of $X_q/\langle\tau\rangle$ is zero.
In particular, $X_q$ is a semi-hyperelliptic curve.
\end{cor}
The next section, we find equations of $X_q$ except for constant numbers by using Theorem \ref{main1} for these $q$.
\section{Equations of $X_q$}
We consider rotation numbers of the automorphism $\tau_n:X_q\ni [z]\mapsto [z+n]\in X_q$.
%
Since $\Gamma_q^n$ acts freely on $\mathbb{H}$ for $q\ge4$, it is sufficient to evaluate rotation numbers at only elements of $S_q$.
\begin{lem}\label{Srot}
Let $n$ be a divisor of $q\ge5$ and $p=\frac{q}{n}$.
We describe by $\frac{x}{z}$ a representative of $\kappa\in S_q$, where $x$ and $z$ are coprime integers.
We take integers $y$ and $w$ such that $xw-yz=1$.
Furthermore, let $k$ be the remainder of $w^2$ divided by $\g(p,z)$.
Then the rotation number at $\kappa$ of $\tau_n$ is
\[{\cal R}_{\tau_n}(\kappa)={\cal R}\left(\frac{p}{\g(p,z)},k\right)\,.\]
\begin{proof}
The width of $\Gamma_q$ is always $q$ since $\Gamma_q$ is a normal subgroup of $SL\left(2,\mathbb{Z}\right)$.
By Lemma \ref{width}, the smallest number $m$ such that $\tau_n^m(\kappa)=\kappa$ satisfied with $p\cdot{\cal W}_{\Gamma_q^n}(\kappa)=m\cdot q$.
By Step 1 of Proposition \ref{h}, we have $m=\frac{p}{\g(p,z)}$.

We recall that elements of $SL\left(2,\mathbb{Z}\right)$ acts continuously on $\hat{\mathbb{H}}$ and it is easy to be calculated the rotation number of parallel displacement at infinity points.
If ${\cal W}_\Gamma(\infty)=R$ for a subgroup $\Gamma\subset SL\left(2,\mathbb{Z}\right)$,
\[\varphi_{[\infty]}(t)=
\begin{cases}
\exp\left(\frac{2\pi it}{R}\right)&(\,{\rm Im}\,t>C\,)\\[1mm]
0&(\,t=\infty\,)
\end{cases}\]
is a chart around the infinity point $[\infty]\in\hat{\mathbb{H}}/\Gamma$.
We take $\left(\begin{smallmatrix}
w&-y\\-z&x
\end{smallmatrix}\right) \in SL\left(2,\mathbb{Z}\right)$ 
which maps $\frac{x}{z}$ to $\infty$.
By easy computation, we have 
\begin{align*}
\begin{pmatrix}w&-y\\-z&x\end{pmatrix}
\left(\,\frac{x}{z}+\varepsilon\,\right)
&=
\frac{1+\varepsilon zw}{-\varepsilon z^2}\\[2mm]
\begin{pmatrix}w&-y\\-z&x\end{pmatrix}
\left(\,\frac{x}{z}+\varepsilon+mn\,\right)
&=
\frac{1+\varepsilon zw+mnzw}{-\varepsilon z^2-mnz^2}
\end{align*}
and
\begin{eqnarray}\label{eq:w^2}
\begin{pmatrix}mnzw+1&mnw^2\\-mnz^2&-mnzw+1\end{pmatrix}
\left(\,\frac{\,1+\varepsilon zw\,}{-\varepsilon z^2}\,\right)
=\frac{1+\varepsilon zw+mnzw}{-\varepsilon z^2-mnz^2}\,.
\end{eqnarray}
Since $mnz=\frac{qz}{\g(p,z)}$ is a multiple of $q$, the matrix of (\ref{eq:w^2}) is equal to
\begin{eqnarray*}
\begin{pmatrix}1&mnw^2\\0&1\end{pmatrix}
\end{eqnarray*}
in modular $q$.
Then the rotation number is given by the remainder of $mnw^2\div mn=w^2$ divided by $p\div m=\g(p,z)$.
\end{proof}
\end{lem}

For $q\le10$ or $q=12$, everything to find an equation of $X_q$ except for values of constant numbers is ready now.
We first consider $X_8$.

Since the denominators of the representatives of the branched points of the natural projection $X_8\rightarrow X_8^1$ is not coprime to 8, 
the branched points of it are $[\infty],\,\left[\frac{3}{8}\right],\,\left[\frac{1}{4}\right],\,\left[\frac{3}{4}\right],\,\left[\frac{1}{2}\right],\,\left[\frac{3}{2}\right],\,\left[\frac{5}{2}\right]$ and $\left[\frac{7}{2}\right]$.
Here $\left[\frac{1}{4}\right],\,\left[\frac{3}{4}\right]$ are $\tau_1$ equivalent.
$\left[\frac{1}{2}\right],\,\left[\frac{3}{2}\right],\,\left[\frac{5}{2}\right],\,\left[\frac{7}{2}\right]$ are also $\tau_1$ equivalent.
By Theorem \ref{main1} and Lemma \ref{Srot}, we get table 2 and
\[y^8=(x-a_1)(x-a_2)(x-a_3)^2(x-a_4)^4\]
which gives an equation of $X_8$.
Here, $n,k$ are rotation numbers, and $m$ is an exponent of an equation about $X_8$.
By normalizing, we have
\begin{eqnarray}\label{eq:a}
y^8=x^2(x-1)(x-a).
\end{eqnarray}
\begin{table}[htb]
\begin{center}
\begin{tabular}{c||c|rccc|lc}
$\frac{x}{z}$&$n$&$xw-yz=1$&example of $w^2$&$\g(8,z)$&$k$&\hspace{9mm}conditions of $m$&$m$\\[1mm] \hline \hline
$\frac{1}{0}$&1&$w=1$&1&8&1&$1=\g(8,m)$ and $m\equiv1\,\mbox{mod}\,8$&1\\[1mm]
$\frac{3}{8}$&1&$3w-8z=1$&9&8&1&$1=\g(8,m)$ and $m\equiv1\,\mbox{mod}\,8$&1\\[1mm] \hline
$\frac{1}{4}$&2&$w-4z=1$&1&4&1&$2=\g(8,m)$ and $\frac{m}{2}\equiv1\,\mbox{mod}\,4$&2\\[1mm] \hline
$\frac{1}{2}$&4&$w-2y=1$&1&2&1&$4=\g(8,m)$ and $\frac{m}{4}\equiv1\,\mbox{mod}\,2$&4
\end{tabular}
\caption{Rotation numbers and exponents about $X_8$}
\end{center}
\end{table}
\begin{theo}\label{proX}
An equation of $X_8$ is given by
\begin{eqnarray}\label{1/2}
y^8=x^2(x-1)(x+1).
\end{eqnarray}
\end{theo}
We give the proof of Theorem \ref{proX} in \S4.
\begin{rem}\label{X_8}
The compact Riemann surface defined by
\begin{eqnarray}\label{4}
y^4=x(x-1)(x+1)(x^2+1)^2
\end{eqnarray}
is isomorphic to the compact Riemann surface defined by {\rm (\ref{1/2})}.
An isomorphism is given by
\begin{align*}
\left\{y^8=x^2(x-1)(x+1)\right\}&\longrightarrow\left\{y^4=x(x-1)(x+1)(x^2+1)^2\right\}\\
(x,y)\hspace{14mm}&\longmapsto\hspace{4mm}\left(\frac{\zeta^4y^4}{x(x+1)},\frac{\sqrt[4]{8}y}{\zeta (x+1)}\right)\\
\left(-\frac{x^2-1}{x^2+1},\frac{\sqrt[4]{2}\zeta y}{x^2+1}\right)\hspace{4mm}&\ \reflectbox{$\longmapsto$}\hspace{14mm}(x,y)\hspace{5mm},
\end{align*}
where $\zeta=\zeta_{16}$.
Therefore, the equation {\rm (\ref{4})} also gives $X_8$.
The form $y^4=f(x)$ corresponds to $g_8^2=0$.
\end{rem}
\begin{rem}
By table 3, table 4 and table 5, we also have equations of $X_9,\,X_{10}$ and $X_{12}$.
They are given by
\begin{align*}
y^9&=x(x-1)^3(x-p_1)^3(x-p_2)^4\\
y^{10}&=x(x-1)^2(x-q_1)^5(x-q_2)^5(x-q_3)^8\\
y^{12}&=x(x-1)^2(x-r_1)^3(x-r_2)^3(x-r_3)^4(x-r_4)^4(x-r_5)^6.
\end{align*}
\tables
\cite{Ya} gives these equations completely.
We see that his equation of $X_9$ is different from our one.
These equations are given by 
\begin{align*}
y^6&=x(x^3+1)y^3+x^5(x^3+1)^2\\
y^{10}&=x(x+1)^2(x-1)^8(x^2+x-1)^5\\
y^{12}&=x(x-1)^2(x+1)^6(x^2+1)^4(x^2-x+1)^3.
\end{align*}
\end{rem}
Of course, the sums of the column of $m$ in these tables are multiple of each $q$ since the infinity points are non-branched points before normalization.
We also obtain an equation of $X_q$ for $q\le7$ by the same way.
In particular, it means that it gives other way to get the classical equation $y^7=x(x-1)^2$ of the Klein's quartic $X_7$.
\section{A canonical model of $X_8$ in the projective space}
As we announced, in this section, we prove Theorem \ref{proX}.
Namely, we determine a constant number $a$ of (\ref{eq:a}).
By Proposition \ref{proX}, two points $(1,0)$ and $(a,0)$ on this algebraic curve are corresponding to $[\infty]$ and $\left[\frac{3}{8}\right]$ on $X_8=\hat{\mathbb{H}}/\Gamma_8$, respectively.
Since there is an automorphism which $[\infty]$ maps $\left[\frac{3}{8}\right]$, for example 
$\left(\begin{smallmatrix}
3&1\\8&3
\end{smallmatrix}\right)$, 
we take an automorphism $\sigma$ which satisfies $\sigma\left((1,0)\right)=(a,0)$.
However, depending the value of $a$, the compact Riemann surface defined by (\ref{eq:a}) dosen't always have such automorphisms.
We see that there is such $\sigma$ if and only if $a=-1$ and thus, we determine $a$ as $-1$.

To find automorphisms of (\ref{eq:a}), we consider a canonical model in the projective space because it is difficult to find automorphisms remains of two variables irreducible polynomials.
The next lemma is fundamental and useful to look for the automorphisms (cf. \cite{KK}).
\begin{lem}
Let $X$ be a non-hyperelliptic compact Riemann surface of genus $g\ge3$ and $X'$ be a canonical model of $X$ in the projective space $\mathbb{P}^{g-1}$.
Then an automorphism $\sigma$ of $X$ is obtained as projective transformation of $\mathbb{P}^{g-1}$ restricted to $X'$.
\end{lem}
We see that $X_8$ is a non-hyperelliptic curve.
If $X_8$ is a hyperelliptic curve, there is an automorphism with order 2 which lies in the center of ${\rm Aut}(X_8)$.
By appendix, ${\rm Aut}(X_8)$ is isomorphic to $PSL\left(2,\mathbb{Z}/8\mathbb{Z}\right)$.
However, the center of $PSL\left(2,\mathbb{Z}/8\mathbb{Z}\right)$ is trivial.
It is a contradiction.

Since $g_8=5$, the projective space is $\mathbb{P}^4$.
We should find a basis of holomorphic differentials of $X_8$ to get a canonical model.
We set projections ${\bf x}:(x,y)\mapsto x$ and ${\bf y}:(x,y)\mapsto y$.
By \S1, we get table 6, here $l=1,2$ and $l'=1,2,3,4$, and a basis as
\[\left<\,\frac{1}{{\bf y}^3}d{\bf x},\frac{{\bf x}}{{\bf y}^5}d{\bf x},\frac{{\bf x}}{{\bf y}^6}d{\bf x},\frac{{\bf x}({\bf x}-1)}{{\bf y}^7}d{\bf x},\frac{{\bf x}}{{\bf y}^7}d{\bf x}\right>.\]
\begin{table}[htb]
\begin{center}
\begin{tabular}{c||c|c|c|c||c|c|c|c|c}
$$&${\bf x}$&${\bf x}-1$&${\bf y}$&$d{\bf x}$
&$\frac{1}{{\bf y}^3}d{\bf x}$&$\frac{{\bf x}}{{\bf y}^5}d{\bf x}$&$\frac{{\bf x}}{{\bf y}^6}d{\bf x}$&$\frac{{\bf x}({\bf x}-1)}{{\bf y}^7}d{\bf x}$&$\frac{{\bf x}}{{\bf y}^7}d{\bf x}$\\[1mm] \hline \hline
$0_l$&4&0&1&3&0&2&1&0&0\\[1mm] \hline
$(1,0)$&0&8&1&7&4&2&1&8&0\\[1mm] \hline
$(a,0)$&0&0&1&7&4&2&1&0&0\\[1mm] \hline
$\infty_{l'}$&$-2$&$-2$&$-1$&$-3$&0&0&1&0&2
\end{tabular}
\caption{Orders of meromorphic functions and differentials}
\end{center}
\end{table}

By three equations
\begin{align*}
\left(\frac{x}{y^6}\right)^2&=\frac{x}{y^5}\cdot\frac{x}{y^7}\\
\left(\frac{x}{y^5}\right)^2&=\frac{1}{y^3}\left(\frac{x(x-1)}{y^7}+\frac{x}{y^7}\right)\\
\left(\frac{1}{y^3}\right)^2&=\frac{x(x-1)}{y^7}\left(\frac{x(x-1)}{y^7}-(a-1)\cdot\frac{x}{y^7}\right),
\end{align*}
we get a canonical model
\begin{eqnarray}\label{pro8}
\left\{[z_1,z_2,z_3,z_4,z_5]\in\mathbb{P}^4\,|\,z_3^2=z_2z_5,z_2^2=z_1(z_4+z_5),z_1^2=z_4\left(z_4-(a-1)z_5\right)\right\}.
\end{eqnarray}
An isomorphism from the algebraic curve (\ref{eq:a}) is
\begin{align*}
(x,y)&\mapsto\left[\frac{1}{y^3},\frac{x}{y^5},\frac{x}{y^6},\frac{x(x-1)}{y^7},\frac{x}{y^7}\right]\\
0_l&\mapsto[\pm\sqrt{a},0,0,-1,1]\\
(1,0)&\mapsto[0,0,0,0,1]\\
(a,0)&\mapsto[0,0,0,a-1,1]\\
\infty_{l'}&\mapsto[1,\pm1,0,1,0],[1,\pm\sqrt{-1},0,-1,0].
\end{align*}
We define an automorphism $\sigma:[z_1,\cdots,z_5]\mapsto[z_1',\cdots,z_5']$ of (\ref{pro8}) by
\[\begin{pmatrix}
z_1'\\\vdots\\z_5'
\end{pmatrix}=
\begin{pmatrix}
c_{1,1}&\cdots&c_{1,5}\\
\vdots&\ddots&\vdots\\
c_{5,1}&\cdots&c_{5,5}
\end{pmatrix}
\begin{pmatrix}
z_1\\\vdots\\z_5
\end{pmatrix}
\]
and then we must have
\begin{align}
z_3'^2&=z_2'z_5'\label{eq:z1}\\
z_2'^2&=z_1'(z_4'+z_5')\label{eq:z2}\\
z_1'^2&=z_4'\left(z_4'-(a-1)z_5'\right).\label{eq:z3}
\end{align}
We should consider the case of an automorphism $\sigma$ maps $[0,0,0,0,1]$ to $[0,0,0,a-1,1]$ and so
\[\lambda\begin{pmatrix}
0\\0\\0\\a-1\\1
\end{pmatrix}=
\begin{pmatrix}
c_{1,1}&\cdots&c_{1,5}\\
\vdots&\ddots&\vdots\\
c_{5,1}&\cdots&c_{5,5}
\end{pmatrix}
\begin{pmatrix}
0\\0\\0\\0\\1
\end{pmatrix}.
\]
Here, $\lambda$ is a non-zero constant.
We may assume $\lambda=1$ and then we have
\[c_{1,5}=c_{2,5}=c_{3,5}=0,\hspace{10pt}
c_{4,5}=a-1,\hspace{10pt}
c_{5,5}=1.\]
\begin{lem}
If an automorphism $\tau$ of $X_8$ satisfies $\tau\left([\infty]\right)=\left[\frac{3}{8}\right]$, we have
\begin{align*}
\tau\left(\left[\frac{3}{8}\right]\right)&=[\infty],\\
\tau\left(\left\{\left[\frac{1}{4}\right],\left[\frac{3}{4}\right]\right\}\right)&=\left\{\left[\frac{1}{4}\right],\left[\frac{3}{4}\right]\right\},\\
\tau\left(\left\{\left[\frac{1}{2}\right],\left[\frac{3}{2}\right],\left[\frac{5}{2}\right],\left[\frac{7}{2}\right]\right\}\right)&=\left\{\left[\frac{1}{2}\right],\left[\frac{3}{2}\right],\left[\frac{5}{2}\right],\left[\frac{7}{2}\right]\right\}.
\end{align*}
\begin{proof}
We regard $\tau$ as a element of $PSL\left(2,\mathbb{Z}/8\mathbb{Z}\right)$.
Since $\tau\left([\infty]\right)=\left[\frac{3}{8}\right]$, $\tau$ is the form
$\pm\left(\begin{smallmatrix}
3&*\\0&3
\end{smallmatrix}\right)$ 
in modular 8.
Then we have this claim by direct calculation.
\end{proof}
\end{lem}
Since $\sigma$ maps $[0,0,0,a-1,1]$ to $[0,0,0,0,1]$ and $a\ne1$, we have
\[c_{1,4}=c_{2,4}=c_{3,4}=0,\hspace{10pt}c_{4,4}=-1.\]
Since $\left[\frac{1}{2}\right]$ corresponding to $\infty_{l'}$ and so $[1,\pm1,0,1,0]$ or $[1,\pm\sqrt{-1},0,-1,0]$, we get
\begin{align*}
c_{3,1}\pm c_{3,2}&=0\\
c_{5,1}\pm c_{5,2}+c_{5,4}&=c_{5,1}\pm \sqrt{-1}c_{5,2}-c_{5,4}=0.
\end{align*}
and therefore, we have
\[c_{3,1}=c_{3,2}=c_{5,1}=c_{5,2}=c_{5,4}=0.\]
By the behavior of the automorphism at $\left[\frac{1}{4}\right]$ corresponding to $0_l$, we also have $c_{2,1}=0$.
Then (\ref{eq:z1}) become to
\begin{align*}
&\,c_{3,3}^2z_3^2=(c_{2,2}z_2+c_{2,3}z_3)(c_{5,3}z_3+z_5)\\
\Leftrightarrow&\ (c_{3,3}^2-c_{2,2}-c_{2,3}c_{5,3})z_3^2=c_{2,2}c_{5,3}z_2z_3+c_{2,3}z_3z_5
\end{align*}
and thus, we have
\[c_{2,3}=0,\hspace{10pt}c_{3,3}^2=c_{2,2}\ne0,\hspace{10pt}c_{5,3}=0.\]
(\ref{eq:z2}) become to
\[c_{2,2}^2z_2^2=(c_{1,1}z_1+c_{1,2}z_2+c_{1,3}z_3)\left(c_{4,1}z_1+c_{4,2}z_2+c_{4,3}z_3-z_4+az_5\right).\]
By the coefficient of $z_3z_5$ and $a\ne0$, we have $c_{1,3}=0$ and then, by the coefficient of $z_2z_5$, we have $c_{1,2}=0$.
Thus, $c_{1,1}\ne0$ and (\ref{eq:z2}) is
\[(c_{1,1}+c_{2,2}^2)z_1z_4=c_{1,1}c_{4,1}z_1^2+c_{1,1}c_{4,2}z_1z_2+c_{1,1}c_{4,3}z_1z_3+(ac_{1,1}-c_{2,2}^2)z_1z_5.\]
We have
\[c_{4,1}=c_{4,2}=c_{4,3}=0,\hspace{10pt}
c_{1,1}=-c_{2,2}^2,\hspace{10pt}a=-1.\]
Therefore, the proof of Theorem \ref{proX} is completed.

We see the form of $\sigma$.
Since we get $c_{1,1}^2=1$ by (\ref{eq:z3}), the form of $\sigma$ is
\[\begin{pmatrix}
-\eta_1&&&&\\
&\eta_2&&&\\
&&\eta_3&&\\
&&&-1&-2\\
&&&&1
\end{pmatrix}.\]
Here, $\eta_1^2=1,\eta_2^2=\eta_1$ and $\eta_3^2=\eta_2$.
We remark that the number of such $\sigma$ is 8 and it is equal to the number of the automorphisms of $X_8$ which maps $[\infty]$ to $\left[\frac{3}{8}\right]$.
$\sigma$ corresponding to $(x,y)\mapsto\left(-x,\eta_3y\right)$ on the semi-hyperelliptic curve $y^8=x^2(x-1)(x+1)$.
We hope that we get equations of $X_9,X_{10}$ and $X_{12}$ completely by like way as $X_8$.
\appendix
\section{Appendix}
\renewcommand{\theequation}{A.\arabic{equation}}
We note that $PSL(q)$ is $PSL\left(2,\mathbb{Z}/q\mathbb{Z}\right)$.
In this appendix, we see some properties of ${\rm Aut}(X_q)$, especially their orders.
We first show ${\rm Aut}(X_q)$ is isomorphic to $PSL(q)$ for $q\ge7$.
Of course, the condition $q\ge7$ is because of $g_q>1$. 
The number of $PSL(q)$ is $R_q<\infty$.
Since elements of $PSL(q)$ are regarded as elements of ${\rm Aut}(X_q)$, we may show $\#{\rm Aut}(X_q)\le R_q$.
We use Hurwitz theorem.
\begin{theo}\label{hurwitz}
Let $X$ be a compact Riemann surface with genus $g>1$ and $\left\{p_1,\cdots,p_n\right\}$ be a maximal set of fixed points of ${\rm Aut}(X)$ inequivalent under the action of ${\rm Aut}(X)$.
We denote the number of the stabilizer of $p_i$ in ${\rm Aut}(X)$ by $m_i$.
Then we get
\begin{eqnarray}\label{eq:hurwitz}
2g-2=N\left(2\overline{g}-2+\sum_{i=1}^n\left(1-\frac{1}{m_i}\right)\right)
\end{eqnarray}
where we have denoted $N=\#{\rm Aut}(X)$ and $\overline{g}$ is the genus of $X/{\rm Aut}(X)$.
\end{theo}
In our case, by Theorem \ref{fre}, we have $g=g_q=1+\frac{(q-6)q^2}{24}\pp$ and $\overline{g}=0
$.
Since the number of the stabilizer of the infinity point in ${\rm Aut}(X_q)$ is at least $q$, we can assume $m_1\ge q$.
We also have $n=3$.
Indeed, if $n\le 2$, we have $2g_q-2<N(-2+1+1)\Leftrightarrow g_q<1$ by (\ref{eq:hurwitz}).
On the other hand, if $n\ge4$, (\ref{eq:hurwitz}) gives
\begin{align*}
&2g_q-2\ge N\left(-2+(n-1)\cdot\frac{1}{2}+1-\frac{1}{q}\right)>N\left(\frac{1}{6}-\frac{1}{q}\right)\\
\Rightarrow\ & N<\frac{q^3}{2}\pp=R_q\,.
\end{align*}
Then (\ref{eq:hurwitz}) is equivalent to
\begin{eqnarray}\label{eq:hurwitz2}
2g_q-2=N\left(1-\frac{1}{m_1}-\frac{1}{m_2}-\frac{1}{m_3}\right)
\end{eqnarray}
in our case.
We can assume $m_2\ge3$ since if $m_2=m_3=2$, we get $g_q<1$ by (\ref{eq:hurwitz2}). 
Then (\ref{eq:hurwitz2}) become to
\[2g_q-2\ge N\left(1-\frac{1}{q}-\frac{1}{3}-\frac{1}{2}\right)
\Leftrightarrow N\le R_q\]
and so $N=R_q$.
Therefore, ${\rm Aut}(X_q)$ is isomorphic to $PSL(q)$ for $q\ge7$.

Next, we consider the largest order of elements of ${\rm Aut}(X_q)$.
If the remainder of $q$ divided by 4 is 2 and isn't divided by 3, we call $q$ type I.
Otherwise, we call $q$ type I\hspace{-.1em}I.
For example, $q=2,10,14,20,\cdots$ are type I.
\begin{lem}
If $q$ is type I, the largest order of elements of $PSL(q)$ is $\frac{3}{2}q$.
If $q$ is type I\hspace{-.1em}I, the largest order is $q$.
\begin{proof}
Since the order of $\left(\begin{smallmatrix}
1&1\\0&1
\end{smallmatrix}\right)\in  PSL(q)$ is $q$ and the order of $\left(\begin{smallmatrix}
p+1&1\\p&1
\end{smallmatrix}\right)\in  PSL(2p)$ is $3p$ for type I $q=2p$, the existence is shown.
We take $A\in PSL(q)$ and see that its order is at most $\frac{3}{2}q$ or $q$.
We prove for each case of $q$.\\[1mm]
{\it Case 1: $q=p$ is prime.}\\
For $p=2$, since $PSL(2)$ is the dihedral group $D_3$, the largest order is 3.
We set $p$ be an odd prime.
$\mathbb{Z}/p\mathbb{Z}$ is a finite field $\mathbb{F}_p$.
Let $\alpha\in\overline{\mathbb{F}}_p$ be an eigenvalue of $A$, where $\overline{\mathbb{F}}_p$ is an algebraic closure of $\mathbb{F}_p$.
$\alpha$ is a solution of
\begin{eqnarray}\label{eq:eigen}
x^2-{\rm tr}(A)x+1=0.
\end{eqnarray}
If $A$ is a diagonalization impossible, $\alpha$ is a multiple root of (\ref{eq:eigen}) and thus, Jordan normal form of $A$ is $\left(\begin{smallmatrix}
1&1\\0&1
\end{smallmatrix}\right)$ whose order is $q$.
We assume $A$ is a diagonalizable matrix.
By Frobenius endomorphism, we realize that $\alpha^p$ is also a solution of (\ref{eq:eigen}).
If $\alpha\ne\alpha^p$, 
$\left(\begin{smallmatrix}
\alpha&0\\0&\alpha^p
\end{smallmatrix}\right)$ is a diagonal matrix of $A$ and $\alpha\cdot\alpha^p=1$.
We have its order is at most $\frac{p+1}{2}$.
If $\alpha=\alpha^p$, the other eigenvalue of $A$ is $\alpha^{-1}$.
Thus, the order of $A$ is at most $\frac{p-1}{2}$.\\[1mm]
{\it Case 2: $q=4$.}\\
By (\ref{eq:eigen}), we have
\begin{align*}
A^3&=\left(t^2-1\right)A-tI\\
A^4&=t\left(t^2-2\right)A-\left(t^2-1\right)I
\end{align*}
where $t={\rm tr}(A)$ and $I$ is the unit matrix.
By considering each case of $t$ in modular 4, we see that the order of $A$ is at most 4.\\[1mm]
{\it Case 3: $q=p^r$ is a prime power.}\\
We use a induction on $r$.
Let assume that the claim is held for $r-1$.
We regard $A\in PSL(p^r)$ as an element of $PSL(p^{r-1})$.
By assumption, we take $n$ such that $A^n$ is a unit matrix in $PSL(p^{r-1})$ with $1\le n\le p^{r-1}$.
It means $A^n=I+p^{r-1}B$ in $PSL(p^r)$ with some $B\in PSL(p^r)$.
Then by
\[A^{np}=\left(I+p^{r-1}B\right)^p=I,\]
we get the order of $A$ is at most $np\left(\le p^r\right)$.\\[1mm]
{\it Case 4: For general $q$.}\\
By Chinese remainder theorem, we may prove it for only $q=2\cdot3^r$, which is type I\hspace{-.1em}I.
If the order of
\[A=B\otimes C\in PSL(2)\otimes PSL(3^r)\]
is larger than $q$, the order of $B$ is 3 and the order of $C$ is larger than $2\cdot3^{r-1}$.
By proof of Case 3, we have that the order of $C$ is $3^r$.
However, then the order of $A=B\otimes C$ is $3^r$.
It is a contradiction.
\end{proof}
\end{lem}
From above results, we have
\begin{prop}\label{most}
For $q\ge7$, the order of elements of ${\rm Aut}(X_q)$ is at most $\frac{3}{2}q$ if $q$ is type I.
If $q$ is type I\hspace{-.1em}I, the order is at most $q$.
\end{prop}
In Theorem \ref{p}, we have the genus of $X_q^1=X_q/\langle z\mapsto z+1\rangle$.
For type I $q=2p$, we shall consider the genus of
\[X_q':=X_q/\left\langle z\mapsto Mz\right\rangle,\]
where
\[M=\begin{pmatrix}p+1&1\\p&1\end{pmatrix},\]
which is an order $3p$ automorphism of $X_q$.
An odd number $p$ isn't divided by 3.
We denote the genus of $X_q'$ by $g_q'$.
\begin{prop}
Let $q=2p$ be type I.
For $q\ge10$, thus for $p\ge5$, we have
\[g_q'=g_{2p}'=1+\frac{\left(p-3{\cal N}\hspace{-1pt}\left(p\right)\right)p}{12}\prod_{l\in{\cal P}(p)}\left(1-\frac{1}{l^2}\right)\,.\]
\begin{proof}
By
\[\begin{pmatrix}p+1&1\\p&1\end{pmatrix}^3
\equiv\begin{pmatrix}1&p+3\\0&1
\end{pmatrix}\hspace{10pt}({\rm mod}\ 2p),\]
$X_q'$ is isomorphic to
\[X_q^2/\left\langle z\mapsto Mz\right\rangle\]
and the order of $M$ is 3 in $X_q^2$.
We see that $X_q^2\rightarrow X_q'$ is the unbranched natural projection.
The subgroup of $SL\left(2,\mathbb{Z}\right)$ generated by $\Gamma_q^2$ and $A$ acts freely on $\mathbb{H}$ since diagonal components of its elements are 1 in modular $p$.
Then we may consider whether only points of $S_q^2$ are branched points or not.
If $\left[\frac{x}{z}\right]\in S_q^2$ with $\g(x,z)=1$ is a branched point of the natural projection, we have
\[\left[\frac{(p+1)x+z}{px+z}\right]_{\Gamma_q^2}=\left[\frac{x}{z}\right]_{\Gamma_q^2}.\]
We remark that $\g\left((p+1)x+z,px+z\right)=\g(x,z)=1$.
Then we have $px+z\equiv z$ or $px+z\equiv-z$ in modular $2p$.
In either case, $x$ is an even number and thus, $z$ is an odd number.
We have
\[\left[\frac{(p+1)x+z}{px+z}\right]_{\Gamma_q^2}=\left[\frac{x+z}{z}\right]_{\Gamma_q^2}\]
and there is $n\in\mathbb{Z}$ such that $x+z\equiv x+2nz$ in modular $2p$.
However, it is a contradiction with $z$ is odd.
 
Since the natural projection $X_q^2\rightarrow X_q'$ is unbranched, we have
\[2g_q^2-2=3\left(2g_q'-2\right)\]
by Hurwitz Theorem \ref{hurwitz}.
The proof is completed by Theorem \ref{p} and a direct computation.
\end{proof}
\end{prop}
By table 7, we notice that $g_q'$ isn't always smaller than $g_q^1$ even though the order of the corresponding automorphism of $X_q'$ is larger than one of $X_q^1$.
Moreover, unfortunately $g_q'\ne0$ except for $q=2$.
Therefore, considering $g_q'$ is of no use to get a semi-hyperelliptic curve after all and we naturally think that $X_q$ is a semi-hyperelliptic curve if and only if $q\le10$ or $q=12$.
\begin{table}[htb]
\begin{center}
\begin{tabular}{c||c|c|c|c|c|c|c|c}
$q$&$2$&10&14&22&26&34&38&$\cdots$\\[1mm] \hline \hline
$g_q$&0&13&49&241&421&1009&1441&$\cdots$\\[1mm] \hline
$g_q^1$&0&0&1&6&10&21&28&$\cdots$\\[1mm] \hline
$g_q'$&0&1&2&6&9&17&22&$\cdots$
\end{tabular}
\caption{Genera for type I $q$}
\end{center}
\end{table}

\small{
}
\end{document}